\documentclass[12pt,oneside,reqno]{amsart}
\usepackage{txfonts}
\usepackage{bbm}
\usepackage{amsmath}
\usepackage{graphicx}
\usepackage{mathrsfs}
\usepackage{stmaryrd}
\usepackage{amsfonts}
\usepackage{enumerate,amsmath,amssymb,amsthm}

\numberwithin{equation}{section}

\newcommand{\be}{\begin{eqnarray}}
\newcommand{\ee}{\end{eqnarray}}
\newcommand{\ce}{\begin{eqnarray*}}
\newcommand{\de}{\end{eqnarray*}}
\newtheorem{theorem}{Theorem}[section]
\newtheorem{lemma}[theorem]{Lemma}
\newtheorem{remark}[theorem]{Remark}
\newtheorem{definition}[theorem]{Definition}
\newtheorem{proposition}[theorem]{Proposition}
\newtheorem{Examples}[theorem]{Example}
\newtheorem{corollary}[theorem]{Corollary}

\def\eps{\varepsilon}

\def\p{\partial}

\def\[{{\Big[}}
\def\]{{\Big]}}
\def\<{{\langle}}
\def\>{{\rangle}}
\def\({{\Big(}}
\def\){{\Big)}}

\def\bx{{\mathbf{x}}}

\def\dif{{\mathord{{\rm d}}}}

\def\no{\nonumber}
\def\={&\!\!=\!\!&}
\def\bt{\begin{theorem}}
\def\et{\end{theorem}}
\def\bl{\begin{lemma}}
\def\el{\end{lemma}}
\def\br{\begin{remark}}
\def\er{\end{remark}}

\def\bd{\begin{definition}}
\def\ed{\end{definition}}
\def\bp{\begin{proposition}}
\def\ep{\end{proposition}}
\def\bc{\begin{corollary}}
\def\ec{\end{corollary}}
\def\bx{\begin{Examples}}
\def\ex{\end{Examples}}

\def\cL{{\mathcal L}}

\def\mE{{\mathbb E}}

\def\mH{{\mathbb H}}

\def\mN{{\mathbb N}}

\def\mQ{{\mathbb Q}}
\def\mR{{\mathbb R}}
\def\mS{{\mathbb S}}

\def\mW{{\mathbb W}}

\def\sB{{\mathscr B}}

\def\sF{{\mathscr F}}

\def\geq{\geqslant}
\def\leq{\leqslant}

\allowdisplaybreaks

\begin{document}

\title{Derivative formula and gradient estimate for SDEs
driven by $\alpha$-stable processes$^*$}

\date{}
\author{Xicheng Zhang}

%\thanks{{\it Keywords: }Reisz transform, critical parabolic equation}

\thanks{$*$This work is supported by NSFs of China (No. 10971076) and Program
for New Century Excellent Talents in University (NCET 10-0654).}

\address{
School of Mathematics and Statistics,
Wuhan University, Wuhan, Hubei 430072, P.R.China\\
Email: XichengZhang@gmail.com
 }

\begin{abstract}
In this paper we prove a derivative formula of Bismut-Elworthy-Li's type as well as gradient estimate
for stochastic differential equations driven by
$\alpha$-stable noises, where $\alpha\in(0,2)$. As an application, the strong Feller property
for stochastic partial differential equations driven by subordinated cylindrical Brownian motions
is presented.
\end{abstract}

\maketitle \rm

\section{Introduction and Main Result}

The derivative formula of diffusion semigroups about the stochastic differential equations
(SDEs) on Riemannian  manifolds was originally introduced by Bismut in \cite{Bi},
and used to derive the estimates of heat kernels and large deviation principles.
His approach is based upon the Malliavin calculus. In \cite{El-Li}, Elworthy and Li used
simple martingale arguments to derive a derivative formula for a large class of diffusion semigroups
on Riemannian manifolds. Nowadays, this type of formula has been proved to be a quite useful tool in various aspects such as
functional inequalities, heat kernel estimates, strong Feller properties and sensitivity analysis, see
\cite{Ar-Th-Wa}, \cite{Ba-Ba-Me}, \cite{Fo-La-Le-Li-To}, \cite{Gu-Wa}, \cite{Pe-Za}, \cite{Wa-Zh} etc.

Due to various successful applications, recently, there is great interest to establish an analogous
derivative formula for jump-diffusion processes. In \cite{Ca-Fr}, Cass and Fritz
proved a derivative formula of Bismut-Elworthy-Li's type for SDEs with jumps and nondegenerate Brownian diffusion term.
In  \cite{Ta}, Takeuchi proved a similar formula for some pure-jump diffusions with finite moments of all orders.
However, their works rule out the interesting $\alpha$-stable processes.
In \cite{Wa1}, Wang used the coupling method and Girsanov's transform to
prove a derivative formula for linear SDEs
driven by pure-jump L\'evy processes including $\alpha$-stable processes,
where the explicit gradient estimates and heat kernel inequalities are also derived.
For infinitely dimensional Ornstein-Uhlenbeck processes with cylindrical $\alpha$-stable processes,
using the infinitely dimensional analysis, Priola and Zabczyk \cite{Pr-Za} established
an explicit derivative formula in terms of the distributional density of $\alpha$-stable processes.
Then the strong Feller property was obtained for a class of semilinear stochastic partial differential equations (SPDEs)
driven by cylindrical $\alpha$-stable processes, where $\alpha\in(1,2)$.

In this paper, we aim to establish a derivative formula of Bismut-Elworthy-Li's type for {\it nonlinear}
SDEs driven by $\alpha$-subordinated Brownian motions. Before moving on, we first recall the classical derivative formulas
for diffusion semigroups. Let $\{W_t\}_{t\geq 0}$ be a standard  $d$-dimensional
Wiener process. Consider the following SDE in $\mR^d$:
\begin{align}
\dif X_t(x)=b_t(X_t(x))\dif t+\sigma\cdot W_t,\ \ X_0(x)=x,\label{St0}
\end{align}
where  $\sigma$ is a $d\times d$ invertible matrix,
and $b:[0,\infty)\times\mR^d\to\mR^d$ satisfies that
\begin{enumerate}[{\bf (H)}]
\item  $b$ has continuous first order partial derivatives with respect to $x$, and
$$
\|\nabla b\|_\infty<+\infty,
$$
where $\nabla b_s(x):=(\p_{x^1} b_s(x),\cdots,\p_{x^d}b_s(x))$ and $\|\cdot\|_\infty$ denotes the uniform norm with
respect to $s$ and $x$.
\end{enumerate}
It is well-known that there are at least two forms of derivative formulas:
for any function $f\in C^1_b(\mR^d)$ and $h\in\mR^d$ (see \cite{Gu-Wa, Zh1}),
\begin{align}
\nabla_h\mE f(X_t(x))=\frac{1}{t}\mE\left(f(X_t(x))\int^t_0\sigma^{-1}
[h+(t-s)\nabla_hb_s(X_s(x))]\dif W_s\right)\label{For1}
\end{align}
and (see \cite{Da-Za, El-Li})
\begin{align}
\nabla_h\mE f(X_t(x))=\frac{1}{t}\mE\left(f(X_t(x))\int^t_0\sigma^{-1}\cdot\nabla_h X_s(x)\dif W_s\right),\label{For2}
\end{align}
where for a function $\varphi$, $\nabla_h\varphi:=\<\nabla\varphi,h\>$ denotes the directional derivative
along $h$, and $\nabla_h X_t(x)$ satisfies the following linear equation:
\begin{align}
\nabla_h X_t(x)=h+\int^t_0\nabla b_s(X_s(x))\cdot\nabla_h X_s(x)\dif s.\label{EE1}
\end{align}
The difference between (\ref{For1}) and (\ref{For2}) lies in that
in formula (\ref{For1}), $\nabla b$ is allowed to be polynomial growth
as in \cite{Gu-Wa, Wa-Zh, Zh1}, while in formula (\ref{For2}), $\nabla b$ usually needs to be bounded.
Below, we shall see that formula (\ref{For2}) is crucial for us since $b$ does not appear explicitly in (\ref{For2}).

Now we turn to the case of SDEs driven by $\alpha$-stable processes.
For $\alpha\in(0,2)$, let $\{S_t\}_{t\geq 0}$ be an independent $\alpha/2$-stable subordinator, i.e., an increasing
$\mR$-valued process with stationary independent increments, and
$$
\mE \mathrm{e}^{\mathrm{i}uS_t}=\mathrm{e}^{t|u|^{\alpha/2}}.
$$
It is well-known that the subordinated Brownian motion $\{W_{S_t}\}_{t\geq 0}$ is an $\alpha$-stable process
(cf. \cite{Be, Sa}).
Let us consider the following SDE in $\mR^d$ driven by $W_{S_t}$:
\begin{align}
\dif X_t(x)=b_t(X_t(x))\dif t+\sigma\cdot \dif W_{S_t},\ \ X_0(x)=x.\label{St}
\end{align}
The generator of Markov process $\{X_t(x),t\geq 0,x\in\mR^d\}$ is given by
$$
\cL_t f(x):=\mathrm{P.V.}\int_{\mR^d} (f(x+\sigma y)-f(x))\frac{\dif y}{|y|^{d+\alpha}}+b_t(x)\cdot\nabla f(x),
$$
where $\mathrm{P.V.}$ stands for the Cauchy principal value. Notice that if $b_t(x)=b_0(x)$ is time homogenous,
then $u_t(x):=\mE f(X_t(x))$ solves the following integro-differential equation (cf. \cite{Ca-Fr, Zh2}):
$$
\p_t u_t(x)=\cL_0 u(x),\ \ u_0=f.
$$
The main aim of the paper is to derive a formula for $\nabla\mE f(X_t(x))$ as stated follows:
\bt\label{Main}
Under {\bf (H)}, for any function $f\in C^1_b(\mR^d)$ and $h\in\mR^d$, we have
\begin{align}
\nabla_h\mE f(X_t(x))=\mE\left(\frac{1}{S_t}f(X_t(x))
\int^t_0\sigma^{-1}\cdot\nabla_h X_s(x)\dif W_{S_s}\right),\label{For}
\end{align}
where $\nabla_h X_s(x)$ is determined by equation (\ref{EE1}). In particular, for any $\alpha\in(0,2)$ and
$p\in(1,\infty]$, there exists a constant $C=C(\alpha,p)>0$ such that for all $t>0$,
\begin{align}
|\nabla\mE f(X_t(x))|\leq C\|\sigma^{-1}\|\mathrm{e}^{\|\nabla b\|_\infty t}t^{-\frac{1}{\alpha}}
\Big(\mE |f(X_t(x))|^p\Big)^{\frac{1}{p}},\label{Gr}
\end{align}
where $\|\sigma^{-1}\|:=\sup_{|x|=1}|\sigma^{-1} x|$ and $|\cdot|$ denotes the Euclidian norm.
\et
\br
From equation (\ref{EE1}), it is easy to see that $s\mapsto\nabla_h X_s(x)$ is a bounded and
continuous $\sigma\{W_{S_r}:r\leq s\}$-adapted process. Thus, the stochastic integral in (\ref{For}) makes sense.
\er

We now introduce the main idea of proving this theorem.
Let $\mW$ be the space of all continuous functions from $[0,\infty)$ to $\mR^d$ vanishing at starting point $0$,
which is endowed with the locally uniform convergence topology and the Wiener measure $\mu_\mW$
so that the coordinate process
$$
W_t(w)=w_t
$$
is a standard $d$-dimensional Brownian motion.
Let $\mS$ be the space of all increasing and c\`adl\`ag functions from $(0,\infty)$ to $(0,\infty)$ with
$\lim_{s\downarrow 0}\ell_s=0,$ which is endowed with the Skorohod metric and
the probability measure $\mu_\mS$ so that the coordinate process
$$
S_t(\ell):=\ell_t
$$
is an $\alpha/2$-stable subordinator $S_t$ (cf. \cite{Be,Sa}).
Consider the following product probability space
$$
(\Omega,\sF,P):=(\mW\times \mS, \sB(\mW)\times\sB(\mS), \mu_\mW\times\mu_\mS),
$$
and define
$$
L_t(w,\ell):=w_{\ell_t}.
$$
Then $\{L_t\}_{t\geq 0}$ is an $\alpha$-stable process on $(\Omega,\sF,P)$. We shall use the following two natural
filtrations associated to the L\'evy process $L_t$ and the Brownian motion $W_t$:
$$
\sF_t:=\sigma\{L_s(w,\ell): s\leq t\},\ \ \sF^\mW_t:=\sigma\{W_s(w): s\leq t\}.
$$
In particular, we can regard the solution $X_t(x)$ of SDE (\ref{St}) as an ($\sF_t$)-adapted
functional on $\Omega$, and therefore,
$$
\mE f(X_t(x))=\int_\mS\int_\mW f(X_t(x; w_{\ell_\cdot}))\mu_\mW(\dif w)\mu_\mS(\dif\ell).
$$
For $\ell\in\mS$, let $X^\ell_t(x)$ solve the following SDE:
\begin{align}
\dif X^\ell_t(x)=b_t(X^\ell_t(x))\dif t+\sigma\cdot \dif W_{\ell_t},\ \ X^\ell_0(x)=x.\label{SDE}
\end{align}
Now, our task is to establish a formula for $\nabla_h\mE^{\mu_\mW} f(X^\ell_t(x))$.
This is not obvious since $t\mapsto W_{\ell_t}$ is not continuous and the classical Bismut-Elworthy-Li formula
can not be used directly.

The remainder of this paper is organized as follows: In Section 2, we shall prove a formula for $\nabla_h\mE^{\mu_\mW} f(X^\ell_t(x))$
by suitable approximation for $\ell_t$. In Section 3, we prove Theorem \ref{Main}.
In Section 4, we prove the strong Feller property for nonlinear
SPDEs driven by subordinated cylindrical Brownian motions.

\section{Derivative formula of SDEs under nonrandom time changed}

In this section, we fix an $\ell\in\mS$ and consider SDE (\ref{SDE}).
If there is no special declaration, all expectations are taken on the Wiener space $(\mW,\sB(\mW),\mu_\mW)$.
First of all, notice that $t\mapsto W_{\ell_t}$ is a Gaussian process with zero means and independent increments.
In particular,
$$
W_{\ell_t}\mbox{ is a c\`adl\`ag $\sF^\mW_{\ell_t}$-martingale.}
$$
Thus,  under {\bf (H)}, it is well-known that for each $x\in\mR^d$,
SDE (\ref{SDE}) admits a unique c\`adl\`ag ($\sF^\mW_{\ell_t}$)-adapted
solution $X^\ell_t(x)$ (cf. \cite[p.249, Theorem 6]{Pr}).

The main aim of this section is to establish the following formula:
\bt\label{Th2}
Under {\bf (H)}, for any  function $f\in C^1_b(\mR^d)$ and $h\in\mR^d$,
we have
\begin{align}
\nabla_h\mE  f(X^\ell_t(x))=\mE \left(\frac{1}{\ell_t}f(X^\ell_t(x))
\int^t_0\sigma^{-1}\cdot\nabla_h X^\ell_s(x)\dif W_{\ell_s}\right),\label{For3}
\end{align}
where $\nabla_h X^\ell_s(x)$ is determined by the following linear equation:
\begin{align}
\nabla_h X^\ell_t(x)=h+\int^t_0\nabla b_s(X^\ell_s(x))\cdot\nabla_h X^\ell_s(x)\dif s.\label{EU222}
\end{align}
\et

For proving this formula, we shall use the time changed argument to transform SDE (\ref{SDE})
into an SDE driven by standard Brownian motions, and then use the classical Bismut-Elworthy-Li formula (\ref{For2}).
For this aim, for $\eps\in(0,1)$, we define
\begin{align}
\ell^\eps_t:=\frac{1}{\eps}\int^{t+\eps}_t\ell_s\dif s+\eps t=\int^1_0\ell_{\eps s+t}\dif s+\eps t.\label{Def}
\end{align}
Since $t\mapsto\ell_t$ is increasing and right continuous, it follows that for each $t\geq 0$,
\begin{align}
\ell^\eps_t\downarrow\ell_t\ \ \mbox{ as \ \  $\eps\downarrow 0$}.\label{ER1}
\end{align}
Moreover, $t\mapsto\ell^\eps_t$ is absolutely continuous and strictly increasing.
Let $\gamma^\eps$ be the inverse function of $\ell^\eps$, i.e.,
$$
\ell^\eps_{\gamma^\eps_t}=t,\ t\geq\ell^\eps_0\ \ \mbox{ and }\ \ \gamma^\eps_{\ell^\eps_t}=t, \ t\geq 0.
$$
By definition, $\gamma^\eps_t$ is also absolutely continuous on
$[\ell^\eps_0,\infty)$.
Let us now define
$$
Y^{\ell^\eps}_t(x):=X^{\ell^\eps}_{\gamma^\eps_t}(x),\ \  t\geq \ell^\eps_0.
$$
By equation (\ref{SDE}) and the change of variables, one sees that for $t\geq \ell^\eps_0$,
\begin{align*}
Y^{\ell^\eps}_t(x)&=x+\int^{\gamma^\eps_t}_0b_s(X^{\ell^\eps}_s(x))\dif s+\sigma\cdot W_t\\
&=x+\int^t_{\ell^\eps_0}b_{\gamma^\eps_s}(Y^{\ell^\eps}_s(x))\dot\gamma^\eps_s\dif s+\sigma\cdot W_t.
\end{align*}
Hence, one can use the classical Bismut-Elworthy-Li formula (see (\ref{For2})) to derive that
$$
\nabla_h\mE  f(Y^{\ell^\eps}_t(x))=\frac{1}{t}\mE\left(f(Y^{\ell^\eps}_t(x))
\int^t_0\sigma^{-1}\cdot\nabla_hY^{\ell^\eps}_s(x)\dif W_s\right), \ t\geq \ell^\eps_0,
$$
where $\nabla_hY^{\ell^\eps}_t(x)$ satisfies
$$
\nabla_hY^{\ell^\eps}_t(x)=h+\int^t_{\ell^\eps_0}\nabla b_{\gamma^\eps_s}(Y^{\ell^\eps}_s(x))
\nabla_hY^{\ell^\eps}_s(x)\dot\gamma^\eps_s\dif s.
$$
Clearly, for each $t\geq 0$,
$$
Y^{\ell^\eps}_{\ell^\eps_t}(x)=X^{\ell^\eps}_t(x),\ \ \
\nabla_hY^{\ell^\eps}_{\ell^\eps_t}(x)=\nabla_hX^{\ell^\eps}_t(x),
$$
and therefore,
\begin{align}
\nabla_h\mE  f(X^{\ell^\eps}_t(x))&=\frac{1}{\ell^\eps_t}\mE
\left(f(X^{\ell^\eps}_t(x))\int^{\ell^\eps_t}_0\sigma^{-1}\cdot
\nabla_hY^{\ell^\eps}_s(x)\dif W_s\right)\no\\
&=\frac{1}{\ell^\eps_t}\mE \left(f(X^{\ell^\eps}_t(x))
\int^t_0\sigma^{-1}\cdot\nabla_hX^{\ell^\eps}_s(x)\dif W_{\ell^\eps_s}\right).\label{Eq1}
\end{align}

Now we want to take limits for both sides of the above formula. We need several lemmas. First of all,
the following lemma is easy.
\bl
For any $p\geq 1$ and $t\geq 0$, we have
\begin{align}
\lim_{\eps\downarrow 0}\mE \left(\sup_{x\in\mR^d}|X^{\ell^\eps}_t(x)-X^\ell_t(x)|^p\right)=0,\label{ER2}
\end{align}
and for any $x, h\in\mR^d$,
\begin{align}
\lim_{\eps\downarrow 0}\mE\left(\sup_{s\in[0,t]}|\nabla_hX^{\ell^\eps}_s(x)-\nabla_hX^\ell_s(x)|^p\right)=0,\label{ER3}
\end{align}
where
\begin{align}
\nabla_h X^{\ell^\eps}_t(x)=h+\int^t_0\nabla b_s(X^{\ell^\eps}_s(x))\cdot\nabla_h X^{\ell^\eps}_s(x)\dif s.\label{ER4}
\end{align}
\el
\begin{proof} For simplicity of notation, we drop the variable ``$x$'' below.
From equation (\ref{SDE}), we have
$$
|X^{\ell^\eps}_t-X^\ell_t|\leq \|\nabla b\|_\infty\int^t_0|X^{\ell^\eps}_s-X^\ell_s|\dif s
+|\sigma\cdot W_{\ell^\eps_t}-\sigma\cdot W_{\ell_t}|.
$$
By Gronwall's inequality, we get
$$
|X^{\ell^\eps}_t-X^\ell_t|\leq \mathrm{e}^{\|\nabla b\|_\infty t}|\sigma\cdot W_{\ell^\eps_t}-\sigma\cdot W_{\ell_t}|,
$$
which then gives (\ref{ER2}) by (\ref{ER1}).

As for (\ref{ER3}), by (\ref{ER4}) and (\ref{EU222}) we have
\begin{align*}
|\nabla_hX^{\ell^\eps}_t-\nabla_hX^\ell_t|
&\leq \int^t_0|\nabla b_s(X^{\ell^\eps}_s)|\cdot|\nabla_hX^{\ell^\eps}_s-\nabla_hX^\ell_s|\dif s
+\int^t_0|\nabla b_s(X^{\ell^\eps}_s)-\nabla b_s(X^\ell_s)|\cdot|\nabla_hX^\ell_s|\dif s,
\end{align*}
which yields by Gronwall's inequality that
$$
|\nabla_hX^{\ell^\eps}_t-\nabla_hX^\ell_t|
\leq \mathrm{e}^{\|\nabla b\|_\infty t}\int^t_0|\nabla b_s(X^{\ell^\eps}_s)-\nabla b_s(X^\ell_s)|\cdot|\nabla_hX^\ell_s|\dif s.
$$
Moreover, from equation (\ref{EU222}), it is easy to see that for any $p\geq 1$,
$$
\sup_{s\in[0,t]}\mE |\nabla_hX^\ell_s|^p\leq C.
$$
Limit (\ref{ER3}) now follows by the dominated convergence theorem, (\ref{ER2})  and the continuity of $x\mapsto\nabla b_s(x)$.
\end{proof}

We also need the following lemma.
\bl\label{Le5}
(i) Assume that $\xi_t$ is a bounded continuous and ($\sF^\mW_{\ell_t}$)-adapted $\mR^d$-valued process.
For each $p,T>0$, we have
\begin{align}
\lim_{\eps\downarrow 0}\mE\left|\int^T_0\xi_s\dif W_{\ell^\eps_s}
-\int^T_0\xi_s\dif W_{\ell_s}\right|^p=0,\label{EU1}
\end{align}
where $\ell^\eps_s$ is defined by (\ref{Def}).

(ii)  Assume that $\xi_t$ is a left continuous and ($\sF^\mW_{\ell_t}$)-adapted $\mR^d$-valued process
and satisfies that for some $p>0$,
\begin{align}
\mE\left(\int^T_0|\xi_s|^2\dif \ell_s\right)^{\frac{p}{2}}<+\infty,\ \ \forall T\geq 0.\label{ET8}
\end{align}
Then there exists a constant $C_p>0$ such that for all $T\geq 0$,
\begin{align}
\mE\left(\sup_{t\in[0,T]}\left|\int^t_0\xi_s\dif W_{\ell_s}\right|^p\right)\leq
C_p\mE\left(\int^T_0|\xi_s|^2\dif \ell_s\right)^{\frac{p}{2}}.\label{EU2}
\end{align}
\el
\begin{proof}
Without loss of generality, we assume $T=1$ and $d=1$.

(i) For $n\in\mN$, set $t_k:=k/n$, $k=0, 1, \cdots, n$ and define
$$
\xi^n_s=\sum_{k=0}^{n-1}\xi_{t_k}1_{[t_k,t_{k+1}]}(s).
$$
By the continuity of $t\mapsto\xi_t$, it is clear that
\begin{align}
\lim_{n\to\infty}\sup_{s\in[0,1]}|\xi^n_s-\xi_s|=0.\label{Cov}
\end{align}
Set
\begin{align*}
\eta_{n,\eps}:=\int^1_0\xi^n_s\dif W_{\ell^\eps_s},&\qquad  \eta_\eps:=\int^1_0\xi_s\dif W_{\ell^\eps_s},\\
\eta_n:=\int^1_0\xi^n_s\dif W_{\ell_s},&\qquad  \eta:=\int^1_0\xi_s\dif W_{\ell_s}.
\end{align*}
Noticing that $W_{\ell^\eps_t}$ is a continuous ($\sF^\mW_{\ell^\eps_t}$)-martingale and
$$
[W_{\ell^\eps}]_t=\ell^\eps_t,
$$
by Burkholder's inequality (cf. \cite[p. 279, Proposition 15.7]{Ka}) and (\ref{ER1}), we have for any $p>0$,
\begin{align*}
\mE|\eta_{n,\eps}-\eta_\eps|^p
=\mE\left|\int^1_0(\xi^n_s-\xi_s)\dif W_{\ell^\eps_s}\right|^p
\leq C\mE\left(\int^1_0|\xi^n_s-\xi_s|^2\dif \ell^\eps_s\right)^{\frac{p}{2}}
\leq C(\ell^1_1)^{\frac{p}{2}}\mE\left(\sup_{s\in[0,1]}|\xi^n_s-\xi_s|^p\right).
\end{align*}
Hence, by the dominated convergence theorem and (\ref{Cov}),
$$
\lim_{n\to\infty}\sup_{\eps\in(0,1)}\mE|\eta_{n,\eps}-\eta_\eps|^p=0.
$$
On the other hand, for fixed $n\in\mN$, by (\ref{ER1}) we have
\begin{align*}
\mE|\eta_{n,\eps}-\eta_n|^p
&=\mE\left|\sum_{k=0}^{n-1}\xi_{t_k}(W_{\ell^\eps_{t_{k+1}}}-W_{\ell^\eps_{t_k}}
-W_{\ell_{t_{k+1}}}+W_{\ell_{t_k}})\right|^p\\
&\leq C_{n,p}K^p\sum_{k=0}^{n-1} \mE|W_{\ell^\eps_{t_{k+1}}}-W_{\ell^\eps_{t_k}}
-W_{\ell_{t_{k+1}}}+W_{\ell_{t_k}}|^p\stackrel{\eps\downarrow 0}{\longrightarrow} 0.
\end{align*}
Combining the above calculations, we obtain (\ref{EU1}).

(ii) We first assume that $\xi_t$ is continuous and bounded by $K$.
Let $\mQ$ be the set of all rational numbers in $[0,1]$. By (\ref{EU1}),
there exists a subsequence $\eps_k\downarrow 0$ such that
$$
\int^t_0\xi_s\dif W_{\ell^{\eps_k}_s}\to \int^t_0\xi_s\dif W_{\ell_s},\ \forall t\in\mQ,\ \ P-a.s.
$$
Therefore,
\begin{align*}
\mE\left(\sup_{t\in[0,1]}\left|\int^t_0\xi_s\dif W_{\ell_s}\right|^p\right)
=\mE\left(\sup_{t\in \mQ}\left|\int^t_0\xi_s\dif W_{\ell_s}\right|^p\right)
=\mE\left(\sup_{t\in \mQ}\lim_{\eps_k\downarrow 0}\left|\int^t_0\xi_s\dif W_{\ell^{\eps_k}_s}\right|^p\right).
\end{align*}
Since $W_{\ell^{\eps_k}_t}$ is a continuous ($\sF^\mW_{\ell^{\eps_k}_t}$)-martingale,
by Fatou's lemma and Burkholder's inequality again, we have
\begin{align}
\mE\left(\sup_{t\in[0,1]}\left|\int^t_0\xi_s\dif W_{\ell_s}\right|^p\right)
&\leq\varliminf_{\eps_k\downarrow0}\mE\left(\sup_{t\in \mQ}\left|\int^t_0\xi_s\dif W_{\ell^{\eps_k}_s}\right|^p\right)\no\\
&\leq C_p\varliminf_{\eps_k\downarrow0}\mE\left(\int^1_0|\xi_s|^2\dif \ell^{\eps_k}_s\right)^{\frac{p}{2}}
=C_p\mE\left(\int^1_0|\xi_s|^2\dif \ell_s\right)^{\frac{p}{2}},\label{EU3}
\end{align}
where the last limit can be proved as in (i).

Next, assume that $\xi_t$ is left-continuous and bounded by $K$. Define
$$
\xi^\eps_t:=\frac{1}{\eps}\int^t_{(t-\eps)\vee 0}\xi_s\dif s,\  \eps>0.
$$
Then $\xi^\eps_t$ is a continuous and bounded ($\sF^\mW_{\ell_t}$)-adapted process. By (\ref{EU3}) we have
$$
\lim_{\eps,\eps'\downarrow 0}\mE\left(\sup_{t\in[0,1]}\left|\int^t_0\xi^\eps_s\dif W_{\ell_s}-
\int^t_0\xi^{\eps'}_s\dif W_{\ell_s}\right|^p\right)
\leq C\lim_{\eps,\eps'\downarrow 0}\left(\int^1_0|\xi^\eps_s-\xi^{\eps'}_s|^2\dif \ell_s\right)^{\frac{p}{2}}=0.
$$
So, (\ref{EU2}) holds for left-continuous and bounded ($\sF^\mW_{\ell_t}$)-adapted process.

In the general case, we truncate $\xi_t$ as follows: for $K>0$, define
$$
\xi^K_t:=(-K)\vee\xi_t\wedge K.
$$
By the dominated convergence theorem, we have
$$
\lim_{K,K'\to\infty}\mE\left(\sup_{t\in[0,1]}\left|\int^t_0\xi^K_s\dif W_{\ell_s}
-\int^t_0\xi^{K'}_s\dif W_{\ell_s}\right|^p\right)
\leq C_p\lim_{K,K'\to\infty}\mE\left(\int^1_0|\xi^K_s-\xi^{K'}_s|^2\dif \ell_s\right)^{\frac{p}{2}}=0.
$$
The proof is complete.
\end{proof}

\br
For $p\geq 1$, by Burkholder's inequality (cf. \cite[p.443, Theorem 23.12]{Ka}), we have
\begin{align}
\mE\left(\sup_{t\in[0,T]}\left|\int^t_0\xi_s\dif W_{\ell_s}\right|^p\right)\asymp
\mE\left(\int^T_0|\xi_s|^2\dif [W_\ell]_s\right)^{\frac{p}{2}},\label{EU22}
\end{align}
where $\asymp$ means that both sides are comparable by multiplying a constant, and
$[W_{\ell}]_t$ denotes the quadratic variation of $(W_{\ell_t})_{t\geq 0}$ given by
$$
[W_{\ell}]_t=\ell_t-\sum_{0<s\leq t}\Delta \ell_s+\sum_{0<s\leq t}|\Delta W_{\ell_s}|^2.
$$
Except for the case of $p=2$, it is not known whether the right hand sides of
(\ref{EU2}) and (\ref{EU22}) are comparable.
\er

\bl
For all $t\geq0$, we have
\begin{align}
\lim_{\eps\downarrow 0}\mE\left(\int^t_0\sigma^{-1}\cdot\nabla_hX^{\ell^\eps}_s(x)\dif W_{\ell^\eps_s}\right)
=\mE\left(\int^t_0\sigma^{-1}\cdot\nabla_hX^\ell_s(x)\dif W_{\ell_s}\right).\label{ER5}
\end{align}
\el
\begin{proof}
For simplicity of notation, we drop the variable ``$x$'' below.
Since $\ell^\eps_s\geq\ell_s$ by definition,
$X^{\ell^\eps}_s$ and $X^\ell_s$ are ($\sF_{\ell^\eps_s}$)-adapted. Thus,
for proving (\ref{ER5}), it suffices to prove the following two limits:
\begin{align}
\lim_{\eps\downarrow 0}\mE\left|\int^t_0(\sigma^{-1}\cdot\nabla_hX^{\ell^\eps}_s
-\sigma^{-1}\cdot\nabla_hX^\ell_s)\dif W_{\ell^\eps_s}\right|^2=0\label{ER6}
\end{align}
and
\begin{align}
\lim_{\eps\downarrow 0}\mE\left|\int^t_0\sigma^{-1}\cdot\nabla_hX^\ell_s\dif W_{\ell^\eps_s}
-\int^t_0\sigma^{-1}\cdot\nabla_hX^\ell_s\dif W_{\ell_s}\right|^2=0.\label{ER7}
\end{align}
For (\ref{ER6}), by the isometry property of stochastic integrals and (\ref{ER3}), we have
\begin{align*}
\lim_{\eps\downarrow 0}\mE\left|\int^t_0\sigma^{-1}\cdot(\nabla_hX^{\ell^\eps}_s
-\nabla_hX^\ell_s)\dif W_{\ell^\eps_s}\right|^2
&=\lim_{\eps\downarrow 0}\mE\left(\int^t_0|\sigma^{-1}\cdot(\nabla_hX^{\ell^\eps}_s-\nabla_hX^\ell_s)|^2\dif\ell^\eps_s\right)\\
&\leq\|\sigma^{-1}\|^2\lim_{\eps\downarrow 0}\mE\left(\sup_{s\in[0,t]}|\nabla_hX^{\ell^\eps}_s-\nabla_hX^\ell_s|^2\right)\ell^1_t=0.
\end{align*}
Limit (\ref{ER7}) follows by Lemma \ref{Le5}.
\end{proof}

{\bf Proof of Theorem \ref{Th2}}: By (\ref{ER1}), (\ref{ER2}) and (\ref{ER5}), the right hand sides of (\ref{Eq1})
converges to the one of (\ref{ER3}). On the other hand, by (\ref{ER2}) and (\ref{ER3}), we have
$$
\nabla_h\mE f(X^{\ell^\eps}_t(x))=\mE\Big(\nabla f(X^{\ell^\eps}_t(x))\cdot\nabla_h X^{\ell^\eps}_t(x)\Big)
\to \mE\Big(\nabla f(X^\ell_t(x))\cdot\nabla_h X^\ell_t(x)\Big)=\nabla_h\mE f(X^\ell_t(x)).
$$

\section{Proof of Theorem \ref{Main}}

The following lemma follows by the monotone class theorem.
\bl\label{Le1}
Let $t\geq 0$ and $A\in\sF_t$. For any $\ell\in\mS$, we have
$$
\{w\in\mW: w_\ell\in A\}\in\sF^\mW_{\ell_t}.
$$
\el
Below we recall a result due to Gin\'e and Marcus \cite{Gi-Ma}.
\bt
Let $\xi_t$ be a left continuous ($\sF_t$)-adapted $\mR^d$-valued process and satisfies
$$
\int^T_0\mE|\xi_s|^\alpha\dif s<+\infty,\ \ \forall T> 0.
$$
Then  there exists a constant $C=C(\alpha)>0$ such that for all $\lambda>0$ and $T>0$,
$$
P\left\{\sup_{t\in[0,T]}\left|\int^t_0\xi_s\dif W_{S_s}\right|\geq \lambda\right\}
\leq C \lambda^{-\alpha}\int^T_0\mE|\xi_s|^\alpha\dif s.
$$
In particular, for any $p\in(0,\alpha)$ and some $C=C(\alpha,p)>0$,
\begin{align}
\mE\left(\sup_{t\in[0,T]}\left|\int^t_0\xi_s\dif W_{S_s}\right|^p\right)
\leq C\left(\int^T_0\mE|\xi_s|^\alpha\dif s\right)^{\frac{p}{\alpha}}.\label{ET4}
\end{align}
\et
\begin{proof}
Define
$$
\zeta:=\sup_{t\in[0,T]}\left|\int^t_0\xi_s\dif W_{S_s}\right|.
$$
Then for any $\eta\geq 0$,
\begin{align*}
\mE\zeta^p&=p\int^\infty_0\lambda^{p-1}P(\zeta>\lambda)\dif\lambda
=p\left(\int^\eta_0+\int^\infty_\eta\right)\lambda^{p-1}P(\zeta>\lambda)\dif\lambda\\
&\leq p\int^\eta_0\lambda^{p-1}\dif\lambda+pC\int^T_0\mE|\xi_s|^\alpha\dif s
\int^\infty_\eta\lambda^{p-\alpha-1}\dif\lambda\\
&\leq \eta^p+\frac{C p}{p-\alpha}\eta^{p-\alpha}\int^T_0\mE|\xi_s|^\alpha\dif s.
\end{align*}
Taking $\eta=\left(\int^T_0\mE|\xi_s|^\alpha\dif s\right)^{\frac{1}{\alpha}}$, we obtain (\ref{ET4}).
\end{proof}

Below we prove a substitution formula about stochastic integrals with respect
to the subordinated Brownian motion $W_{S_t}$.
\bp
Assume that $\xi_t(W_S)$ is a bounded and left continuous ($\sF_t$)-adapted $\mR^d$-valued process.
Then for any $T\geq 0$, we have
\begin{align}
\int^T_0\xi_s(W_S)\dif W_{S_s}=\int^T_0\xi_s(W_\ell)\dif W_{\ell_s}\Big|_{\ell=S},\ \ P-a.s.\label{ET7}
\end{align}
Moreover, for any nonnegative random variable $g$ on $\mS$ and $p>0$, we have
\begin{align}
\mE\left(g(S)\sup_{t\in[0,T]}\left|\int^t_0\xi_s(W_S)\dif W_{S_s}\right|^p\right)\leq C_p
\int_{\mS}g(\ell)\mE^{\mu_\mW}\left(\int^T_0|\xi_s(W_\ell)|^2\dif \ell_s\right)^{\frac{p}{2}}\mu_\mS(\dif\ell).\label{ET77}
\end{align}
\ep
\begin{proof}
The proof is the same as in Lemma \ref{Le5}. Without loss of generality, we assume $T=1$.
For given $n\in\mN$, set $t_k:=k/n$, $k=0,1,\cdots, n$ and define
$$
\xi^n_s(W_S):=\xi_0(W_S)1_{\{0\}}(s)+\sum_{i=0}^{n-1}\xi_{t_i}(W_S)1_{(t_i, t_{i+1}]}(s).
$$
Then, by definition we have
\begin{align*}
\int^1_0\xi^n_s(W_S)\dif W_{S_s}
&=\sum_{i=0}^{n-1}\xi^n_{t_i}(W_S)(W_{S_{t_{i+1}}}-W_{S_{t_i}})\\
&=\sum_{i=0}^{n-1}\xi^n_{t_i}(W_\ell)(W_{\ell_{t_{i+1}}}-W_{\ell_{t_i}})\Big|_{\ell=S}
=\int^1_0\xi^n_s(W_\ell)\dif W_{\ell_s}\Big|_{\ell=S},
\end{align*}
and by the left continuity of $s\mapsto\xi_s(W_S)$,
\begin{align}
\lim_{n\to\infty}|\xi^n_s(W_S)-\xi_s(W_S)|=0,\ \ s\geq 0.\label{ET5}
\end{align}
For $p\in(0,\alpha)$, by (\ref{ET4}) we have
$$
\mE\left|\int^1_0(\xi^n_s(W_S)-\xi_s(W_S))\dif W_{S_s}\right|^p\leq C\int^1_0\mE|\xi^n_s(W_S)-\xi_s(W_S)|^p\dif s\to 0.
$$
On the other hand, by Lemma \ref{Le1}, for each $\ell\in\mS$ and $s\geq 0$,
$\xi_s(W_\ell)$ is ($\sF^\mW_{\ell_s}$)-adapted. Thus, by (\ref{EU2}) we have
\begin{align*}
\mE\left|\int^1_0(\xi^n_s(W_\ell)-\xi_s(W_\ell))\dif W_{\ell_s}\Big|_{\ell=S}\right|^p
&=\int_\mS\mE^{\mu_{\mW}}\left|\int^1_0(\xi^n_s(W_\ell)-\xi_s(W_\ell))\dif W_{\ell_s}\right|^p\mu_{\mS}(\dif\ell)\\
&\leq C\int_\mS\mE^{\mu_{\mW}}\left(\int^1_0|\xi^n_s(W_\ell)-\xi_s(W_\ell)|^2
\dif \ell_s\right)^{\frac{p}{2}}\mu_{\mS}(\dif\ell)\\
&=C\mE\left(\int^1_0|\xi^n_s(W_S)-\xi_s(W_S)|^2\dif S_s\right)^{\frac{p}{2}}.
\end{align*}
which converges to zero by (\ref{ET5}) and the dominated convergence theorem.
Combining the above calculations, we obtain (\ref{ET7}). Lastly, (\ref{ET77}) is an easy consequence of (\ref{ET7})
and (\ref{EU2}).
\end{proof}

We are now in a position to give

{\bf Proof of Theorem \ref{Main}}:
Formula (\ref{For}) follows by (\ref{For3}) and (\ref{ET7}).
We now prove gradient estimate (\ref{Gr}). By H\"older's inequality, we have
\begin{align*}
|\nabla_h\mE f(X_t(x))|&\leq\mE\left(\frac{1}{S_t}\left|f(X_t(x))\int^t_0\sigma^{-1}\cdot\nabla_h X_s(x)\dif W_{S_s}\right|\right)\\
&\leq\Big(\mE |f(X_t(x))|^p\Big)^{\frac{1}{p}}
\left(\mE\left|\frac{1}{S_t}\int^t_0\sigma^{-1}\cdot\nabla_h X_s(x)\dif W_{S_s}\right|^q\right)^{\frac{1}{q}},
\end{align*}
where $\frac{1}{p}+\frac{1}{q}=1$.
By (\ref{EE1}), it is easy to see that
$$
|\nabla_h X_t(x)|\leq |h| \mathrm{e}^{\|\nabla b\|_\infty t}.
$$
Thus, by (\ref{ET77}) we have
\begin{align*}
\mE\left|\frac{1}{S_t}\int^t_0\sigma^{-1}\cdot\nabla_h X_s(x)\dif W_{S_s}\right|^q
&\leq C_q\int_{\mS}\frac{1}{\ell^q_t}\mE^{\mu_\mW}\left(\int^t_0|\sigma^{-1}
\cdot\nabla_h X_s(x)|^2\dif \ell_s\right)^{\frac{q}{2}}\mu_\mS(\dif\ell)\\
&\leq C_q\|\sigma^{-1}\|^q|h|^q\mathrm{e}^{q\|\nabla b\|_\infty t}\int_{\mS}\frac{\mu_\mS(\dif\ell)}{\ell^{q/2}_t}.
\end{align*}
Recalling that the distributional density of $\alpha$-stable
subordinator satisfies (cf. \cite[(14)]{Bo-St-Sz})
$$
P\circ S^{-1}_t(\dif s)\leq C ts^{-1-\frac{\alpha}{2}}\mathrm{e}^{-ts^{-\frac{\alpha}{2}}}\dif s,
$$
we have
\begin{align*}
\int_{\mS}\frac{\mu_\mS(\dif\ell)}{\ell^{q/2}_t}=\mE\left(\frac{1}{S^{q/2}_t}\right)
\leq C\int^\infty_0ts^{-1-\frac{\alpha+q}{2}}\mathrm{e}^{-ts^{-\frac{\alpha}{2}}}\dif s
=C t^{-\frac{q}{\alpha}}\int^\infty_0 u^{\frac{q}{\alpha}}\mathrm{e}^{-u}\dif u,
\end{align*}
where the last equality is due to the change of variable $u=ts^{-\frac{\alpha}{2}}$, and
$C$ only depends on $\alpha,q$.
Combining the above calculations, we obtain (\ref{Gr}).

\section{Strong Feller property of SPDEs driven by subordinated cylindrical Brownian motions}

Let $\mH$ be a real separable Hilbert space with the inner product $\<\cdot,\cdot\>_\mH$. The norm in $\mH$
is denoted by $\|\cdot\|_\mH$. Let $A$ be a negative self-adjoint operator in $\mH$ with discrete spectrals, i.e.,
there exists an orthogonal basis $\{e_k\}_{k\in\mN}$ and a sequence of real numbers
$0<\lambda_1\leq\lambda_2\leq\cdots\leq\lambda_k\to\infty$ such that
$$
A e_k=-\lambda_k e_k.
$$
Let $\{W^k_t, t\geq 0\}_{k\in\mN}$ be a sequence of independent standard $1$-dimensional Brownian motion.
Let $\{L_t\}_{t\geq 0}$ be the subordinated cylindrical Brownian motion in $\mH$ defined by
$$
L_t:=\sum_{k=1}^\infty\beta_kW^k_{S_t}e_k,\ \ \beta_k\in\mR,
$$
where $S_t$ is an independent $\alpha/2$-stable subordinator.

Consider the following SPDE in Hilbert space $\mH$:
\begin{align}
\dif X_t(x)=[AX_t(x)+F(X_t(x))]\dif t+\dif L_t,\ \ X_0=x\in\mH.\label{Eq}
\end{align}
Our aim of this section is to prove that
\bt\label{Th}
Assume that for some $\delta>0$,
\begin{align}
\beta_k\geq\delta,\ \ \forall k\in\mN\ \ \mbox{ and }\ \
\sum_{k=1}^\infty\frac{\beta^2_k}{\lambda_k}<+\infty,\label{EW1}
\end{align}
and one of the following conditions holds:
\begin{enumerate}[(i)]
\item $\alpha\in(1,2)$ and $F:\mH\to\mH$ is Lipschitz continuous;
\item $\alpha\in(0,2)$ and $F:\mH\to\mH$ is bounded and Lipschitz continuous.
\end{enumerate}
Then there exists a constant $C=C(\alpha)>0$ such that for any bounded Borel measurable
function $f: \mH\to\mR$, $x,y\in\mH$ and $t>0$,
\begin{align}
|\mE f(X_t(x))-\mE f(X_t(y))|\leq C\delta \mathrm{e}^{\|F\|_{\mathrm{Lip}}t}
t^{-\frac{1}{\alpha}}\|f\|_\infty\|x-y\|_\mH,\label{For5}
\end{align}
where $\|F\|_{\mathrm{Lip}}:=\sup_{x\not=y}\frac{\|F(x)-F(y)\|_\mH}{\|x-y\|_\mH}$.
\et

Let us first prove a result about the following stochastic convolution (cf. \cite[Theorem 4.5]{Pr-Za}):
$$
Z^A_t:=\int^t_0\mathrm{e}^{A(t-s)}\dif L_s=\sum_{k=1}^\infty\beta_k\int^t_0 \mathrm{e}^{-\lambda_k(t-s)}\dif W^k_{S_s} e_k.
$$
\bp
If $\sum_{k=1}^\infty\beta^2_k/\lambda_k<+\infty$, then for any $p\in(0,\alpha)$ and $t\geq 0$, we have
\begin{align}
\mE\|Z^A_t\|_\mH^p\leq C_{\alpha,p}\left(\sum_{k=1}^\infty\frac{|\beta_k|^2}
{\lambda_k}\right)^{\frac{\alpha}{2}}t^{1-\frac{\alpha}{2}}.\label{ER8}
\end{align}
\ep
\begin{proof}
Recall the following Gaussian formula (cf. \cite[Appendix A.1]{Nu}): Let $\{\xi_k\}_{k\in\mN}$ be a sequence of independent
random variables defined on some probability space $(\Omega',\sF',P')$ with normal distribution $N(0,1)$,
and $\{c_k\}_{k\in\mN}$ a sequence of real numbers. Then
$$
\mE'\left|\sum_{k=1}^\infty c_k\xi_k\right|^p=A_p
\left(\sum_{k=1}^\infty|c_k|^2\right)^{\frac{p}{2}},\ \ A_p:=\int_\mR\frac{|x|^p}{\sqrt{2\pi}}\mathrm{e}^{-\frac{x^2}{2}}\dif x.
$$
Using this formula and by Fubini's theorem and (\ref{ET4}), we have
\begin{align*}
A_p\mE\|Z^A_t\|^p_\mH&=A_p\mE\left(\sum_{k=1}^\infty\left|\beta_k\int^t_0 \mathrm{e}^{-\lambda_k(t-s)}\dif W^k_{S_s}\right|^2\right)^{\frac{p}{2}}
=\mE\mE'\left|\sum_{k=1}^\infty\xi_k\beta_k\int^t_0 \mathrm{e}^{-\lambda_k(t-s)}\dif W^k_{S_s}\right|^p\\
&=\mE'\mE\left|\sum_{k=1}^\infty\xi_k\beta_k\int^t_0 \mathrm{e}^{-\lambda_k(t-s)}\dif W^k_{S_s}\right|^p\leq
C\mE'\int^t_0\left(\sum_{k=1}^\infty |\xi_k\beta_k|^2\mathrm{e}^{-2\lambda_k(t-s)}\right)^{\frac{\alpha}{2}}\dif s\\
&\leq C\left(\sum_{k=1}^\infty\int^t_0\mE'|\xi_k\beta_k|^2
\mathrm{e}^{-2\lambda_k(t-s)}\dif s\right)^{\frac{\alpha}{2}}t^{1-\frac{\alpha}{2}}
\leq C\left(\sum_{k=1}^\infty\frac{|\beta_k|^2}{2\lambda_k}\right)^{\frac{\alpha}{2}}t^{1-\frac{\alpha}{2}}.
\end{align*}
The proof is finished.
\end{proof}
We now establish the following existence and uniqueness of mild solutions to equation (\ref{Eq}).
\bp
Assume that $F:\mH\to\mH$ is Lipschitz continuous and $\sum_{k=1}^\infty\beta^2_k/\lambda_k<+\infty$.
Then for each $x\in\mH$,  there exists a unique $X_t(x)\in\mH$ satisfying
\begin{align}
X_t=\mathrm{e}^{At}x+\int^t_0\mathrm{e}^{A(t-s)}F(X_s)\dif s+Z^A_t.\label{Eq2}
\end{align}
\ep
\begin{proof}
Consider the following deterministic equation:
$$
Y_t=\mathrm{e}^{At}x+\int^t_0\mathrm{e}^{A(t-s)}F(Y_s+Z^A_s)\dif s.
$$
Using the standard Picard's iteration, it is easy to see that there exists a unique $Y\in C([0,\infty);\mH)$
satisfying the above equation. Thus, $X_t=Y_t+Z^A_t$ satisfies (\ref{Eq2}).
\end{proof}

Let $\mH_n:=\{x=\sum_{k=1}^n c_k e_k, c_k\in\mR\}$ and $\Pi_n$ the projection operator from $\mH$ to $\mH_n$
defined by
$$
\Pi_n x:=\sum_{k=0}^n\<x,e_k\>_\mH e_k.
$$
Let $\rho_n$ be a sequence of nonnegative smooth functions with
$$
\mathrm{supp}(\rho_n)\subset\{z\in\mH_n: |z|\leq 1/n\},\ \ \int_{\mH_n}\rho_n(z)\dif z=1.
$$
Define
$$
F_n(x):=\int_{\mH_n}\rho_n(\Pi_n x-z)\Pi_n F(z)\dif z
$$
and
$$
L^n_t:=\Pi_n L_t=\sum_{k=1}^n\beta_kW^k_{S_t}e_k.
$$
Consider the following finite dimensional approximation of equation (\ref{Eq}):
$$
\dif X^n_t=[\Pi_nAX^n_t+F_n(X^n_t)]\dif t+\dif L^n_t,\ \ X^n_0=\Pi_n x.
$$
\bl
Under the assumptions of Theorem \ref{Th}, for any fixed $t>0$ and $x\in\mH$, we have
\begin{align}
\lim_{n\to\infty}\|X^n_t(\Pi_n x)-X_t(x)\|_\mH=0,\ \ P-a.s.\label{ET9}
\end{align}
\el
\begin{proof}
By Duhamel's formula, one can write
$$
X^n_t=\mathrm{e}^{At}\Pi_n x+\int^t_0\mathrm{e}^{A(t-s)}F_n(X^n_s)\dif s+\Pi_nZ^A_t.
$$
Set
$$
Y^n_t:=X^n_t-\Pi_nZ^A_t,\ \ Y_t:=X_t-Z^A_t.
$$
Then
$$
Y^n_t-Y_t=\mathrm{e}^{At}(\Pi_n x-x)+\int^t_0\mathrm{e}^{A(t-s)}(F_n(Y^n_s+\Pi_nZ^A_s)-F(Y_s+Z^A_s))\dif s.
$$
Hence,
$$
\|Y^n_t-Y_t\|_\mH\leq\|\Pi_n x-x\|_\mH+\int^t_0\|F_n(Y^n_s+\Pi_nZ^A_s)-F(Y_s+Z^A_s)\|_\mH\dif s.
$$
Notice that
\begin{align*}
\|F_n(Y^n_s+\Pi_nZ^A_s)-F(Y_s+Z^A_s)\|_\mH
&\leq \|F\|_{\mathrm{Lip}}(\|Y^n_s-Y_s\|_\mH+\|(\Pi_n-I)Z^A_s\|_\mH)\no\\
&\quad+\|(\Pi_n-I)F(Y_s+Z^A_s)\|_\mH
\end{align*}
and
$$
\lim_{n\to\infty}\|(\Pi_n-I)Z^A_s\|_\mH=0,\ \ \lim_{n\to\infty}\|(\Pi_n-I)F(Y_s+Z^A_s)\|_\mH=0.
$$
(i) If $\alpha\in(1,2)$ and $F$ is Lipschitz continuous,
by (\ref{ER8}) and the dominated convergence theorem, we have
\begin{align}
\varlimsup_{n\to\infty}\|Y^n_t-Y_t\|_\mH\leq C\int^t_0\varlimsup_{n\to\infty}\|Y^n_s-Y_s\|_\mH\dif s,\label{EW2}
\end{align}
which then gives
\begin{align}
\varlimsup_{n\to\infty}\|Y^n_t-Y_t\|_\mH=0.\label{EW3}
\end{align}
(ii) If $\alpha\in(0,2)$ and $F$ is bounded, by Fatou's lemma, we also have (\ref{EW2}) and (\ref{EW3}).
\end{proof}

{\bf Proof of Theorem \ref{Th}}:
For any function $f\in C^1_b(\mH)$, by (\ref{Gr}) with $p=\infty$,
there exists a constant $C=C(\alpha)>0$ such that for all $x,y\in\mH$ and $t>0$,
\begin{align}
|\mE f(X^n_t(\Pi_nx))-\mE f(X^n_t(\Pi_ny))|&\leq C\delta
\|f\|_\infty \mathrm{e}^{\|\nabla F_n\|_\infty t}t^{-\frac{1}{\alpha}}\|\Pi_nx-\Pi_ny\|_\mH\no\\
&\leq C\delta\|f\|_\infty \mathrm{e}^{\|F\|_{\mathrm{Lip}} t} t^{-\frac{1}{\alpha}}\|x-y\|_\mH.\label{For4}
\end{align}
By taking limits for (\ref{For4}), we get (\ref{For5}) for any $f\in C^1_b(\mH)$. For general bounded measurable $f$,
it follows by a standard approximation (see \cite[p.125, Lemma 7.1.5]{Da-Za}).

\vspace{5mm}

{\bf Example:} Consider the following nonlinear stochastic heat equation in $[0,1]$ with Dirichlet boundary conditions:
$$
\p_t u=[\Delta u+ b(u)]\dif t+\dif L_t,\ \ u_t(0)=u_t(1)=0,\ \ u_0=\varphi,
$$
where $b:\mR\to\mR$ is a Lipschitz function, and $\varphi\in L^2([0,1])=:\mH$.
It is well-known that $A=\Delta$ is a negative self-adjoint operator on $\mH$ with eigenvalues
$$
\lambda_k=\pi^2k^2,\ \ k\in\mN,
$$
and eigenvectors
$$
e_k(\zeta)=\sqrt{2}\sin(\pi k\zeta),\ \ \zeta\in[0,1].
$$
In particular,
$$
\sum_{k\in\mN}\frac{1}{\lambda_k}<+\infty.
$$
Thus, if one takes $\beta_k=1$, then (\ref{EW1}) holds.

\vspace{5mm}

{\bf Acknowledgements:}

The author thanks Lihu Xu for useful discussions.


\begin{thebibliography}{999}

\bibitem{Ar-Th-Wa}Arnaudon A., Thalmaier A., Wang F.Y.: Gradient estimate and Harnack inequality
on non-compact Riemannian manifolds. Stoch. Proc. Appl., 119(2009), 3653-3670.

\bibitem{Ba-Ba-Me}Bally V., Bavouzet M.P. and Messaoud M.: Integration by parts formula for locally smooth
laws and applications to sensitivity computations. Ann. Appl. Prob., 17, 33-66(2007).

\bibitem{Be}Beroin J.: L\'evy processes. Cambridge Tracts in Math., Cambridge Univ. Press, London (1996)

\bibitem{Bi}Bismut J.M.: Large deviations and the Malliavin calculus. Birkh\"auser, Boston (1984).

\bibitem{Bo-St-Sz}Bogdan K. and St\'os A. and Sztonyk P.: Harnack inequality for stable processes on $d$-sets.
Studia Math., 158(2), 163-198(2003).

\bibitem{Ca-Fr}Cass T.R. and Friz P.K.: The Bismut-Elworthy-Li formula for jump-diffusions and applications to
Monte-Carlo methods in finance. arXiv: math.PR/0604311 v1.

\bibitem{Da-Za}Da Prato G. and Zabczyk J.: Ergodicity for infinite dimensional systems. London
Math. Society, Lecture Notes Series, 229, Cambridge Univ. Press, 1996.

\bibitem{El-Li}Elworthy K.D. and Li X.M.: Formulae for the derivatives of heat semigroups. J. Func. Anal.,
125, 252-286(1994).

\bibitem{Fo-La-Le-Li-To}Fourni\'e E., Lasry J.M., Lebuchoux J., Lions P.L., Touzi N.: Applications of Malliavin
calculus to Monte-Carlo methods in finance. Finance stoch., 3, 391-412(1999).

\bibitem{Gi-Ma}Gin\'e E. and Marcus B.: The central limit theorem for stochastic integrals with respect to
L\'evy processes. The Annals of Prob., 11, 58-77(1983).

\bibitem{Gu-Wa}Guillin A. and Wang F.Y.: Degenerate Fokker-Planck equations: Bismut formula, gradient estimate and
Harnack inequality. J. Diff. Equa., in press(2012).

\bibitem{Ka}Kallenberg O.: Foundations of Modern Probability. Springer Verlag, 1997.

\bibitem{Nu}Nualart D.: The Malliavin calculus and related topics. Second Edition, Springer-Verlag, Berlin(2005).

\bibitem{Pe-Za}Peszat S. and Zabczyk J.: Strong Feller property and irreducibility for diffusions on Hilbert spaces.
The Annals of Prob., 23, 157-172(1995).

\bibitem{Pr-Za}Priola E. and Zabczyk J.: Structural properties of semilinear SPDEs driven by cylindrical stable
processes. Probab. Theory Related Fields 149 (2011), 97-137.

\bibitem{Pr}Protter P.E.: Stochastic integration and differential equations. Second Edition, Springer-Verlag, Berlin, 2004.

\bibitem{Sa}Sato K.I.: L\'evy processes and infinite divisble distributions. Cambridge University Press, Cambridge, 1999.

\bibitem{Ta}Takeuchi A.: The Bismut-Elworthy-Li type formulae for stochastic differential equations with jumps.
Journal of Theoretical Probability, 23(2010), 576-604.

\bibitem{Wa1}Wang F.Y.: Derivative formula and Harnack inequality for jump processes. http://arxiv.org/abs/1104.5531.

\bibitem{Wa-Zh}Wang F.Y. and Zhang X.C.: Derivative formula and applications for degenerate diffusion semigroups.
http://arxiv.org/abs/1107.0096.

\bibitem{Zh1} Zhang X.C.: Stochastic flows and Bismut formulas for stochastic Hamiltonian systems.
Stoch. Proc. Appl. 120(2010), 1929--1949.

\bibitem{Zh2}Zhang, X.C.: Stochastic functional differential equations driven by L\'evy processes
and quasi-linear partial integro-differential equations. Ann. Appl. Prob. in press(2012).

\end{thebibliography}
\end{document}